\documentclass[11pt]{article}
\usepackage{amsfonts}
\usepackage{mathrsfs}
\usepackage{amsmath}
\usepackage{amssymb}
\usepackage{graphicx}
\usepackage{color}
\usepackage{ifpdf}
\usepackage{epic}
\usepackage{amsthm,latexsym}
\usepackage{indentfirst}
\renewcommand{\paragraph}{\roman{paragraph}}
\def \A{\mathcal{A}}

\def \T{\mathcal{T}}
\def \G{\mathscr{G}}
\usepackage[]{caption2} 
\newtheorem{theorem}{\scshape \mdseries  Theorem}[section]
\newtheorem{lemma}[theorem]{\scshape \mdseries  Lemma}
\newtheorem{coro}[theorem]{\scshape \mdseries  Corollary}

\newtheorem{defi}[theorem]{\scshape \mdseries  Definition}

\begin{document}
\title{\sf Maximizing spectral radii of uniform hypergraphs with few edges\thanks{Supported by 
  National Natural Science Foundation of China (11371028, 71101002), Project of Educational Department of Anhui Province (KJ2012B040),
  Scientific Research Fund for Fostering Distinguished Young Scholars of Anhui University(KJJQ1001),  Academic Innovation Team of Anhui University Project (KJTD001B).}}
\author{Yi-Zheng Fan$^{1,}$\thanks{Corresponding author. E-mail address: fanyz@ahu.edu.cn (Y.-Z. Fan), tansusan1@ahjzu.edu.cn (Y.-Y. Tan)}, Ying-Ying Tan$^{1,2}$, Xi-Xi Peng$^{1}$, An-Hong Liu$^1$ \\
{\small  \it $1$. School of Mathematical Sciences, Anhui University, Hefei 230601, P. R. China}\\
    {\small  \it $2$. School of Mathematics \& Physics, Anhui Jianzhu University, Hefei 230601, P. R. China}
}
\date{}
\maketitle

\noindent {\bf Abstract:} In this paper we investigate the hypergraphs whose spectral radii attain the maximum among all uniform hypergraphs with given number of edges.
In particular we characterize the hypergraph(s) with maximum spectral radius over all unicyclic hypergraphs, linear or power unicyclic hypergraphs with given girth, linear or power bicyclic hypergraphs, respectively.



\section{Introduction}
Let $G$ be a simple graph on $n$ vertices with $m$ edges. The spectral radius of $G$ is the largest eigenvalue of the adjacency matrix of $G$.
In 1985 Brualdi and Hoffman \cite{bh} investigated the maximum spectral radius of the adjacency matrix of
a, not necessarily connected, graph in the set of all graphs with given number of vertices and
edges. Their work was followed by other people, in the connected graph case as well as in the
general case, and a number of papers have been written.
In particular, Rowlinson \cite{row} settled the
problem for the general case; he proved that among all graphs with fixed number of edges (or,
equivalently, with fixed number of vertices and edges), there is a unique graph that maximizes
the spectral radius of the adjacency matrix. (The unique graph turns out to be a threshold graph.)
However, the problem of determining the maximizing graphs, i.e., the connected case
of the problem, is still unresolved, though we know the optimal graph is a maximal graph \cite{ord}, where a graph is called {\it maximal} if its degree sequence is
majorized by no other graphic sequences \cite{merris}.

The  maximizing graphs have been identified only for some choices of $n$ and $m$.
By the definition of maximal graphs, a maximal graph always contains a vertex adjacent to all other vertices.
So, the maximizing graph contains a vertex adjacent to all other vertices, which was proved by Brualdi and Solheid \cite{bs}.
As a conclusion, the maximizing tree of order $n$ is the star $K_{1,n-1}$ \cite{cs} and the maximizing unicyclic graph is obtained from $K_{1,n-1}$ by adding an edge between two pendant vertices \cite{hong}.

In this paper, we consider the similar problem for $k$-uniform hypergraphs, that is,
maximizing the spectral radius of the adjacency tensor of $k$-uniform hypergraphs among all $k$-uniform hypergraphs
with given number of vertices and edges.
A $k$-uniform hypergraph $G=(V,E)$ consists of a vertex set $V$ and an edge set $E \subseteq {V \choose k}$.
The adjacency tensor \cite{cd} of $G$ is defined as the $k$-th order $n$ dimensional tensor $\A(G)$ whose $i_{1}i_{2}\ldots i_{k}$-entry is
$$
  a_{i_{1}i_{2}\ldots i_{k}}=\left\{
  \begin{array}{cl}
   \frac{1}{(k-1)!} & \hbox{if~} \{v_{i_{1}},v_{i_{2}},\ldots,v_{i_{k}}\} \in E(G)\\
  0, &  \hbox{otherwise.}
  \end{array}
    \right.
  $$
Qi \cite{qi05} introduces the eigenvalues of a supersymmetric tensor, from which one can get the definition of the eigenvalues of the adjacency tensor of a $k$-uniform hypergraph.
The spectral radius of a $k$-uniform hypergraph is the maximum modulus of the eigenvalues of its adjacency tensor; see  more in Section 2.

We show that among all connected $k$-uniform hypergraphs with given number of vertices and edges, the one with maximum spectral radius contains a vertex adjacent to all other vertices, which is parallel to the result on simple graphs.
As a conclusion,  among all $k$-uniform hypertrees with given number of edges, the hyperstar is the unique maximizing one.
We determine the maximizing unicyclic $k$-uniform hypergraphs, and characterize the maximizing linear or power unicyclic/bicyclic hypergraphs.
All hypergraphs in this paper are $k$-uniform with $k \ge 3$.

\section{Preliminaries}
Let $G$ be a $k$-uniform hypergraph.
The {\it degree} $d_v$ of a vertex $v \in V(G)$ is defined as $d_v=|\{e_{j}:v\in e_{j}\in E(G)\}|$.
 A {\it walk} $W$ of length $l$ in $G$ is a sequence of alternate vertices and edges: $v_{0}e_{1}v_{1}e_{2}\ldots e_{l}v_{l}$,
    where $\{v_{i},v_{i+1}\}\subseteq e_{i}$ for $i=0,1,\ldots,l-1$.
If $v_0=v_l$, then $W$ is called a {\it circuit}.
A walk of $G$ is called a {\it path} if no vertices or edges are repeated.
A circuit $G$ is called a {\it cycle} if no vertices or edges are repeated except $v_0=v_l$.
The  hypergraph $G$ is said to be {\it connected} if every two vertices are connected by a walk.

If $G$ is connected and acyclic, then $G$ is called a {\it hypertree}.
It is known that a connected $k$-uniform hypergraph with $n$ vertices and $m$ edges is acyclic if and only if $m=\frac{n-1}{k-1}$, i.e. $n=m(k-1)+1$;
see \cite[Proposition 4, p.392]{ber}.
If $G$ is connected and contains exactly one cycle, then $G$ is called a {\it unicyclic hypergraph}.

\begin{lemma}\label{numedges}
If $G$ is a unicyclic $k$-uniform hypergraph with $n$ vertices and $m$ edges, then $n=m(k-1)$.
\end{lemma}

{\bf Proof.}
Let $u_1 e_1 u_2 \cdots u_{t-1} e_t u_1$ be the unique cycle of $G$.
Now adding a new vertex $w$ into $G$, and replacing the edge $e_1$ by $(e_1 \backslash \{u_1\}) \cup \{w\}$,
we will arrive at a new acyclic hypergraph which has $n+1$ vertices and $m$ edges.
So $m=\frac{(n+1)-1}{k-1}=\frac{n}{k-1}$. The result follows. \hfill $\blacksquare$

\begin{defi}\label{cyclo}
Let $G$ be a $k$-uniform hypergraph with $n$ vertices, $m$ edges and $l$ connected components.
The cyclomatic number of $G$ is denoted and defined by $c(G)=m(k-1)-n+l$.
The hypergraph $G$ is called a $c(G)$-cyclic hypergraph.
\end{defi}

If $k=2$, the above definition is exactly that of simple graphs.
In particular, a connected hypergraph $G$ is called {\it bicyclic} if $c(G)=2$.

\begin{figure}
\centering
  \setlength{\unitlength}{1bp}
  \begin{picture}(393.12, 206.35)(0,0)
  \put(0,0){\includegraphics{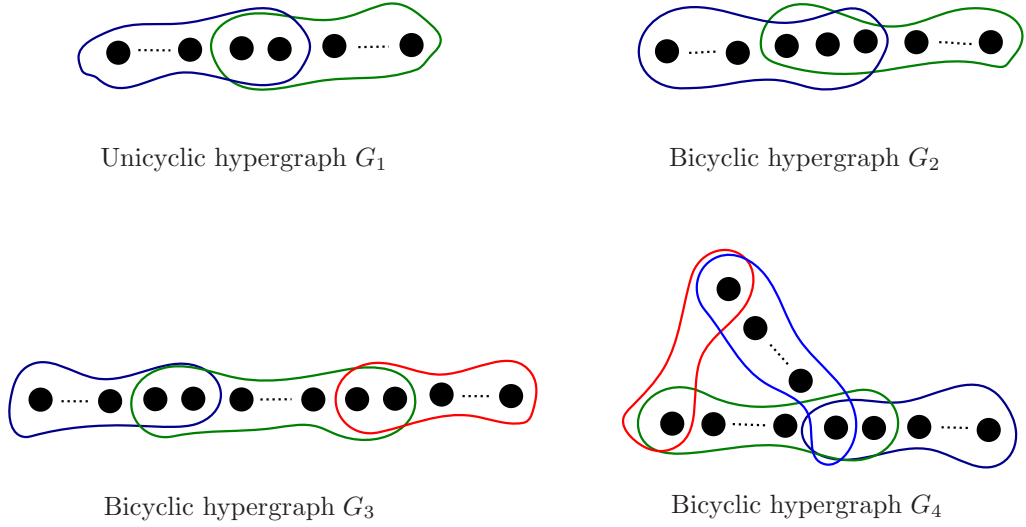}}
  \put(39.95,139.26){\fontsize{10}{12}\selectfont Unicyclic hypergraph $G_1$}
  \put(254.29,139.26){\fontsize{10}{12}\selectfont Bicyclic hypergraph $G_2$}
  \put(41.35,7.87){\fontsize{10}{12}\selectfont Bicyclic hypergraph $G_3$}
  \put(255.17,8.39){\fontsize{10}{12}\selectfont Bicyclic hypergraph $G_4$}
  \end{picture}%
\caption{An illustration of nonlinear unicyclic or bicyclic hypergraphs}\label{uni-bic}
\end{figure}

\begin{figure}
\centering
  \setlength{\unitlength}{1bp}
  \begin{picture}(463.80, 286.53)(0,0)
  \put(0,0){\includegraphics{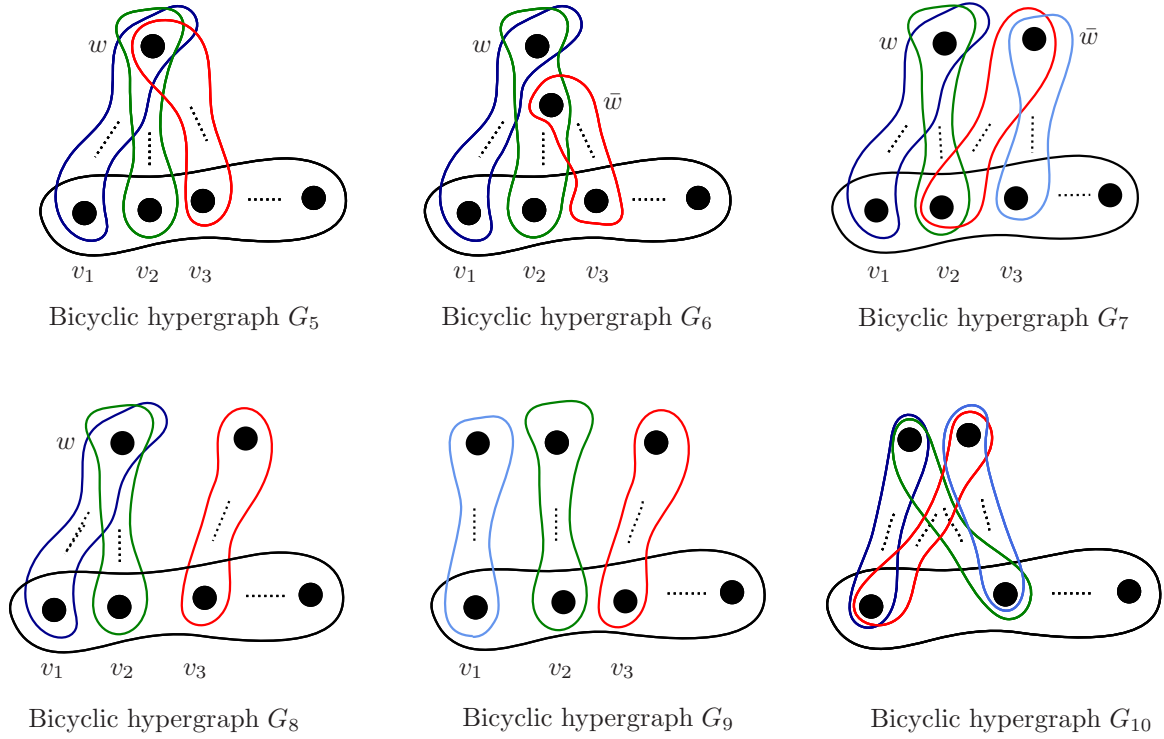}}
  \put(20.25,159.49){\fontsize{10}{12}\selectfont Bicyclic hypergraph $G_5$}
  \put(168.24,159.49){\fontsize{10}{12}\selectfont Bicyclic hypergraph $G_6$}
  \put(329.87,7.93){\fontsize{10}{12}\selectfont Bicyclic hypergraph $G_{10}$}
  \put(12.60,7.87){\fontsize{10}{12}\selectfont Bicyclic hypergraph $G_8$}
  \put(325.28,159.33){\fontsize{10}{12}\selectfont Bicyclic hypergraph $G_7$}
  \put(176.18,8.16){\fontsize{10}{12}\selectfont Bicyclic hypergraph $G_9$}
  \put(35.09,262.39){\fontsize{10}{12}\selectfont $w$}
  \put(179.74,262.39){\fontsize{10}{12}\selectfont $w$}
  \put(23.01,112.68){\fontsize{10}{12}\selectfont $w$}
  \put(332.96,263.56){\fontsize{10}{12}\selectfont $w$}
  \put(229.25,240.95){\fontsize{10}{12}\selectfont $\bar{w}$}
  \put(408.59,265.51){\fontsize{10}{12}\selectfont $\bar{w}$}
  \put(28.86,177.40){\fontsize{10}{12}\selectfont $v_1$}
  \put(53.03,177.40){\fontsize{10}{12}\selectfont $v_2$}
  \put(72.91,177.40){\fontsize{10}{12}\selectfont $v_3$}
  \put(172.72,177.40){\fontsize{10}{12}\selectfont $v_1$}
  \put(198.84,177.40){\fontsize{10}{12}\selectfont $v_2$}
  \put(222.62,177.40){\fontsize{10}{12}\selectfont $v_3$}
  \put(328.67,177.40){\fontsize{10}{12}\selectfont $v_1$}
  \put(354.79,177.40){\fontsize{10}{12}\selectfont $v_2$}
  \put(378.57,177.40){\fontsize{10}{12}\selectfont $v_3$}
  \put(17.16,27.69){\fontsize{10}{12}\selectfont $v_1$}
  \put(43.28,27.69){\fontsize{10}{12}\selectfont $v_2$}
  \put(70.96,27.69){\fontsize{10}{12}\selectfont $v_3$}
  \put(174.67,27.69){\fontsize{10}{12}\selectfont $v_1$}
  \put(208.59,27.69){\fontsize{10}{12}\selectfont $v_2$}
  \put(231.98,27.69){\fontsize{10}{12}\selectfont $v_3$}
  \end{picture}
\caption{An illustration of linear bicyclic hypergraphs}\label{bic-line}
\end{figure}

If $|e_i \cap e_j| \in \{ 0,s\}$ for all edges $e_i \neq e_j$ of a hypergraph $G$, then $G$ is called an {\it $s$-hypergraph}.
A simple graph is a $2$-uniform $1$-hypergraph.
Note that $1$-hypergraphs here are also called {\it linear hypergraphs} \cite{bre}.
So, a hypertree is a linear hypergraph; otherwise, if two edges $e_1,e_2$ have two vertices $v_1,v_2$ in common, then
$v_1e_1v_2e_2v_1$ is a $2$-cycle.
See Fig. \ref{uni-bic} and Fig. \ref{bic-line} for some examples of nonlinear or linear unicyclic/bicyclic uniform hypergraphs.

\begin{defi}{\em \cite{HQS}}
Let $G=(V,E)$ be a simple graph. For any $k\geq3$, the {\it $k$-th power of $G$}, denoted by $G^{k}:=(V^{k},E^{k})$,
  is defined as the $k$-uniform hypergraph with the set of vertices $V^{k}:=V\cup{\{i_{e,1},\ldots,i_{e,k-2}|e\in E}\}$ and the set of edges
$E^{k}:={\{e\cup{{\{i_{e,1},\ldots,i_{e,k-2}}}\}|e\in E}\}$.
\end{defi}

Obviously, the power of simple graphs (or simply called {\it power hypergraphs}) are linear.
A $k$-uniform {\it hyperstar} with $m$ edges is the $k$-th power of the ordinary star $K_{1,m}$,
 and a {\it loose path} with $m$ edges is the $k$-th power of the ordinary path with $m$ edges.

For integers $k\geq 3$ and $n\geq 2$,
  a real {\it tensor} (also called {\it hypermatrix}) $\mathcal{T}=(t_{i_{1}\ldots i_{k}})$ of order $k$ and dimension $n$ refers to a
  multidimensional array with entries $t_{i_{1}i_2\ldots i_{k}}$ such that $t_{i_{1}i_2\ldots i_{k}}\in \mathbb{R}$ for all $i_{j}\in [n]:=\{1,2,\ldots,n\}$ and $j\in [k]$.
 The tensor $\mathcal{T}$ is called \textit{symmetric} if its entries are invariant under any permutation of their indices.
 Given a vector $x\in \mathbb{R}^{n}$, $\mathcal{T}x^{k}$ is a real number, and $\mathcal{T}x^{k-1}$ is an $n$-dimensional vector, which are defined as follows:
   $$\mathcal{T}x^{k}=\sum_{i_1,i_{2},\ldots,i_{k}\in [n]}t_{i_1i_{2}\ldots i_{k}}x_{i_1}x_{i_{2}}\cdots x_{i_k},~
   (\mathcal{T}x^{k-1})_i=\sum_{i_{2},\ldots,i_{k}\in [n]}t_{ii_{2}i_3\ldots i_{k}}x_{i_{2}}x_{i_3}\cdots x_{i_k} \mbox{~for~} i \in [n].$$
 Let $\mathcal{I}$ be the {\it identity tensor} of order $k$ and dimension $n$, that is, $i_{i_{1}i_2 \ldots i_{k}}=1$ if and only if
   $i_{1}=i_2=\cdots=i_{k} \in [n]$ and zero otherwise.

\begin{defi}{\em \cite{qi05}} Let $\mathcal{T}$ be a $k$-th order $n$-dimensional real tensor.
For some $\lambda \in \mathbb{C}$, if the polynomial system $(\lambda \mathcal{I}-\mathcal{T})x^{k-1}=0$, or equivalently $\mathcal{T}x^{k-1}=\lambda x^{[k-1]}$, has a solution $x\in \mathbb{C}^{n}\backslash \{0\}$,
then $\lambda $ is called an eigenvalue of $\mathcal{T}$ and $x$ is an eigenvector of $\mathcal{T}$ associated with $\lambda$,
where $x^{[k-1]}:=(x_1^{k-1}, x_2^{k-1},\ldots,x_n^{k-1}) \in \mathbb{C}^n$.
\end{defi}

If $x$ is a real eigenvector of $\mathcal{T}$, surely the corresponding eigenvalue $\lambda$ is real.
In this case, $x$ is called an {\it $H$-eigenvector} and $\lambda$ is called an {\it $H$-eigenvalue}.
Furthermore, if $x\in \mathbb{R}_{+}^{n}$ (the set of nonnegative vectors of dimension $n$), then $\lambda $ is called an {\it $H^{+}$-eigenvalue} of $\mathcal{T}$;
if $x\in \mathbb{R}_{++}^{n}$ (the set of positive vectors of dimension $n$), then $\lambda$ is said to be an {\it $H^{++}$-eigenvalue} of $\mathcal{T}$.
The {\it spectral radius of $\T$} is defined as
$$\rho(\T)=\max\{|\lambda|: \lambda \mbox{ is an eigenvalue of } \T \}.$$

Chang et al. \cite{CPZ} introduced the irreducibility of tensor. A tensor $\T=(t_{i_{1}\ldots i_{k}})$ of order $k$ and dimension $n$ is called {\it reducible} if there exists a nonempty proper subset $I \subset [n]$ such that
$t_{i_1i_2\ldots i_k}=0$ for any $i_1 \in I$ and any $i_2,\ldots,i_k \notin I$.
If $\T$ is not reducible, then it is called {\it irreducible}.
Friedland et al. \cite{FGH} proposed a weak version of the irreducibility of nonnegative tensors $\T$.
The graph associated with $\T$, denoted by $G(\T)$, is the directed graph with vertices $1, 2, \ldots, n$ and an edge from $i$ to $j$
 if and only if $t_{ii_2\ldots i_k}>0$ for some $i_l = j$, $l = 2, \ldots, m$.
The tensor $\T$ is called {\it weakly irreducible} if $G(\T)$ is strongly connected.
Surely, an irreducible tensor is always weakly irreducible.
Pearson and Zhang \cite{PZ} proved that the adjacency tensor of $\A(G)$ is weakly irreducible if and only if $G$ is connected.

\begin{theorem}\label{PF} {\em \bf (The Perron-Frobenius Theorem for Nonnegative Tensors)}

1. {\em (Yang and Yang 2010 \cite{YY})} If $\T$ is a nonnegative tensor of order $k$ and dimension $n$, then $\rho(\T)$ is an $H^+$-eigenvalue of $\T$.

2. {\em (Frieland, Gaubert and Han 2011 \cite{FGH})} If furthermore $\T$ is weakly irreducible, then $\rho(\T)$ is the unique $H^{++}$-eigenvalue of $\T$,
with the unique eigenvector $x \in \mathbb{R}_{++}^{n}$, up to a positive scaling coefficient.

3. {\em(Chang, Pearson and Zhang 2008 \cite{CPZ})} If moreover $\T$ is irreducible, then $\rho(\T)$ is the unique $H^{+}$-eigenvalue of $\T$,
with the unique eigenvector $x \in \mathbb{R}_{+}^{n}$, up to a positive scaling coefficient.

\end{theorem}

Let $x=(x_{1},x_{2},\ldots,x_{n})^T\in R^{n}$, and let $G$ be a hypergraph on vertices $v_{1},v_{2},\ldots,v_{n}$.
Then $x$ can be considered as a function defined on $V(G)$, that is, each vertex $v_{i}$ is mapped to $x_{i}=:x_{v_{i}}$.
If $x$ is an eigenvector of $\A(G)$, then it defined on $V(G)$ naturally, i.e. $x_v$ is the entry of $x$ corresponding to $v$.
From the Theorem \ref{PF}, the spectral radius of $\A(G)$, also referred to the spectral radius of $G$, denoted by $\rho(G)$, is exactly the
largest $H$-eigenvalue of $\A(G)$.
If $G$ is connected, then there exists a unique positive eigenvector up to scales corresponding to $\rho(G)$, called the {\it Perron vector} of $G$.
In addition, $\rho(G)$ is the optimal value of the following maximization (see \cite{qi05}):
$$ \rho(G)=\max_{x \in \mathbb{R}^n, \|x\|_k=1}\A(G)x^k=\max_{x \in \mathbb{R}^n, \|x\|_k=1}\sum_{e=\{u_1,u_2,\cdots, u_k\} \in E(G)}k x_{u_1}x_{u_2}\cdots x_{u_k}.\eqno(2.1)$$
The eigenvector equation  $\mathcal{A}(G)x^{k-1}=\lambda x^{[k-1]}$ could be interpreted as
$$ \lambda x_u^{k-1}= \sum_{\{u,u_2,u_3,\ldots, u_k\} \in E(G)} x_{u_2}x_{u_3} \cdots x_{u_k}, \mbox{~for each~} u \in V(G).\eqno(2.2)$$

Li, Shao and Qi \cite{lsq} introduce the operation of {\it moving edges} on hypergraphs.
Let $r \ge  1$ and let $G$ be a hypergraph with $u \in V(G)$ and $e_1, \ldots, e_r \in E(G)$ such that $u \notin e_i$ for $i = 1, \ldots, r$.
Suppose that $v_i \in e_i$ and write $e'_i = (e_i \backslash \{v_i\}) \cup \{u\}$ ($i = 1, \ldots, r)$.
Let $G'$ be the hypergraph with $V(G')=V(G)$ and $E(G') = (E \backslash \{e_1, \ldots, e_r\}) \cup \{e'_1, \ldots, e'_r\}$.
We say that $G'$ is obtained from $G$ by moving edges $(e_1, \ldots, e_r)$ from $(v_1, \ldots, v_r)$ to $u$.

\begin{theorem} \label{moveedges} {\em \cite{lsq}}
Let $r \ge  1$ and let $G$ be a connected hypergraph.
 Let $G'$ be obtained from $G$ by moving edges $(e_1, \ldots, e_r)$ from $(v_1, \ldots, v_r)$ to $u$. 
 Assume that $G'$ contains no multiple edges.
If $x$ is a Perron vector of $G$ and $x_u \ge \max_{1\le i \le r}x_{v_i}$, then $\rho(G') > \rho(G)$.
\end{theorem}

\begin{figure}
\centering
  \setlength{\unitlength}{1bp}
  \begin{picture}(411.70, 125.21)(0,0)
  \put(0,0){\includegraphics{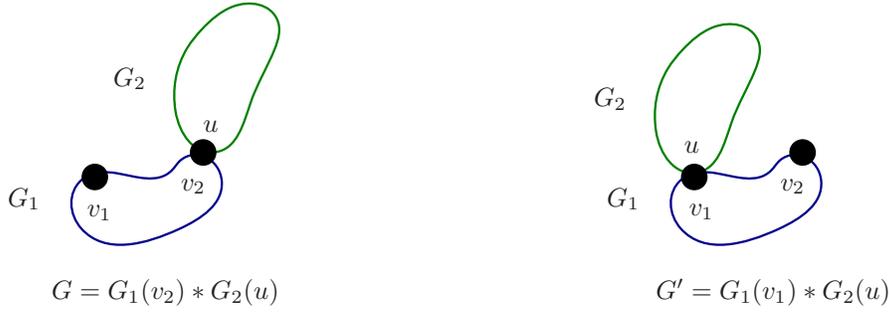}}
  \put(35.58,39.33){\fontsize{10}{12}\selectfont $v_1$}
  \put(70.97,49.11){\fontsize{10}{12}\selectfont $v_2$}
  \put(79.45,71.04){\fontsize{10}{12}\selectfont $u$}
  \put(45.55,88.49){\fontsize{10}{12}\selectfont $G_2$}
  \put(5.67,43.12){\fontsize{10}{12}\selectfont $G_1$}
  \put(262.47,39.63){\fontsize{10}{12}\selectfont $v_1$}
  \put(296.87,49.11){\fontsize{10}{12}\selectfont $v_2$}
  \put(260.57,63.29){\fontsize{10}{12}\selectfont $u$}
  \put(226.68,80.73){\fontsize{10}{12}\selectfont $G_2$}
  \put(231.56,43.12){\fontsize{10}{12}\selectfont $G_1$}
  \put(22.29,8.12){\fontsize{10}{12}\selectfont $G=G_{1}(v_{2})\ast G_{2}(u)$}
  \put(250.00,8.12){\fontsize{10}{12}\selectfont $G'=G_{1}(v_{1})\ast G_{2}(u)$}
  \end{picture}
\caption{An illustration of relocating subhypergraph}\label{mov}
\end{figure}

We now introduce a special case of moving edges.
Let $G_1,G_2$ be two vertex-disjoint hypergraphs, where $v_{1},v_{2}$ are two distinct vertices of $G_{1}$ and $u$ is a vertex of $G_{2}$ (called the root of $G_2$).
Let $G=G_{1}(v_{2})\ast G_{2}(u)$ (respectively, $G'=G_{1}(v_{1})\ast G_{2}(u)$) be the hypergraph obtained by identifying $v_2$ with $u$ (respectively, identifying $v_1$ with $u$);
see the graphs in Fig. \ref{mov}.
We say that $G'$ is obtained from $G$ by {\it relocating} $G_{2}$ rooted at $u$ from $v_{2}$ to $v_{1}$.

\begin{coro}\label{relocate}
Let $G=G_{1}(v_{2})\ast G_{2}(u)$ and $G'=G_{1}(v_{1})\ast G_{2}(u)$ be two connected hypergraphs.
If there exists a Perron vector $x$ of $G$ such that $x_{v_1} \ge x_{v_2}$, then $\rho(G') > \rho(G)$.
\end{coro}

\section{Maximizing the spectral radii of uniform hypergraphs}
If $G$ is a hypergraph whose spectral radius attains the maximum among a certain class of hypergraphs,
then $G$ is called a {\it maximizing hypergraph} in such class.
For a connected acyclic hypergraph (or hypertree), unicyclic or bicyclic $k$-uniform hypergraph, the number of vertices is determined by the number of its edges by Definition \ref{cyclo}.
So we only mention the number of edges of hypertrees, unicyclic or bicyclic hypergraphs in the following discussion.

\begin{lemma} \label{deg}If $G$ is a maximizing hypergraph among the connected hypergraphs with fixed number edges, then $G$ contains a vertex adjacent to all other vertices.
\end{lemma}

\noindent{\bf Proof.} Let $x$ be a Perron vector of $G$.
 By Theorem \ref{PF}, $x$ is positive.
 We take one vertex, say $u_{0}$ of $G$, such that $x_{u_{0}}=\max{\{x_{v}: v\in V(G)}\}$.
 Suppose there exists a vertex $w$ not adjacent to $u_0$.
 As $G$ is connected, there exists a path connecting $u_0$ and $w$, say $u_0 e_1 u_1 \cdots u_{t-1}e_t u_t$, where $t \ge 2$ and $u_t=w$.
 Moving the edge $e_t$ from $u_{t-1}$ to $u_0$, we will arrive at a new hypergraph $G'$ which contains a new edge $e'_t:=(e_t \backslash \{u_{t-1}\}) \cup u_0$.
 Note that $e'_t \notin G$ otherwise $w$ would be adjacent to $u_0$.
 Since $x_{u_0} \ge x_{t-1}$, by Theorem \ref{moveedges}, we get $\rho(G') > \rho(G)$; a contradiction.\hfill $\blacksquare$

By Lemma \ref{deg}, we easily get the following result proved by Li, Shao and Qi \cite{lsq}.

 \begin{coro} {\em \cite{lsq}} Among all hypertrees with $m$ edges, the hyperstar $K_{1,m}^k$ is the unique maximizing hypergraph.
 \end{coro}

 \begin{coro} Among all unicyclic hypergraphs with $m$ edges, the unique maximizing hypergraph is obtained from the hypergraph $G_1$ in Fig. \ref{uni-bic} by attaching a hyperstar $K_{1,m-2}^k$ with its center at one vertex of degree $2$.
 \end{coro}

 \noindent{\bf Proof.}
 Let $G$ be a maximizing hypergraph.
 By Lemma \ref{deg}, $G$ contains a vertex $u_0$ adjacent to all other vertices.
 Let $\bar{G}$ be the sub-hypergraph induced by the edges containing $u_0$. Surely, $V(\bar{G})=V(G)$.
 If $\bar{G}$ is a hyperstar (centered at $u_0$), noting that $G$ is unicyclic, $G$ consists of $\bar{G}$ and an edge only containing vertices of $V(G)\backslash \{u_0\}$.
 Then $G$ would contains more than one cycles; a contradiction.
 So $\bar{G}$ is not a hyperstar, that is, it contains a pair of edges sharing a common vertex except the vertex $u_0$.
The result follows.  \hfill $\blacksquare$

Denote by $\G_m(G_0)$ the class of hypergraphs with $m$ edges each obtained from a fixed connected hypergraph $G_0$ by attaching some hypertrees at some vertices of $G_0$ respectively
(i.e. identifying a vertex of a hypertree with some vertex of $G_0$ each time).
We first discuss the maximizing hypergraph(s) in $\G_m(G_0)$, and then get some corollaries for special hypergraphs.

\begin{lemma} \label{max}
If $G$ is a maximizing connected hypergraph in $\G_m(G_0)$, then $G$ is obtained from $G_0$ by attaching a hyperstar with its center at some vertex $u$ of $G_0$.
Furthermore, if $x$ is a Perron vector of $G$, then $x_u >x_v$ for any other vertex $v$ of $G$;
if $G_0$ contains more than one edge, then $u$ has degree greater than one in $G_0$.
\end{lemma}

\noindent{\bf Proof.}
 Let $x$ be the Perron vector of $G$, and let $u \in V(G_0)$ such that $x_u=\max\{x_v: v \in V(G_0)\}$.
The result will follow by the following two claims.

 {\scshape \mdseries Claim 1}:  All hypertrees are attached at $u$.
 Otherwise, if there exists a hypertrees $T_v$ attached at $v \ne u$ of $G_0$,
relocating $T_v$ from $v$ to $u$, noting that $x_{u} \ge  x_v$, we will get a hypergraph with a larger spectral radius by Corollary \ref{relocate}.
So we assume $G$ is obtained from $G_0$ by attaching exactly one hypertree $T_u$ at $u$.

The above discussion also implies that $u$ is the unique vertex in $G_0$ with maximum value given by $x$.
Furthermore, $u$ is unique vertex of $G$ with maximum value; otherwise, if $\bar{v}$ is one outside $G_0$ such that $x_{\bar{v}} \ge x_u$,
relocating $G_0$ from $u$ to $\bar{v}$,  we also get a  contradiction by Corollary \ref{relocate}.

Suppose that $G_0$ contains more than one edge.
Assume that $u$ has degree $1$ in $G_0$ and lies in some edge $e_0$ of $G_0$.
As $G_0$ is connected, $e_0$ contains a vertex, say $\bar{w}$, with degree at least $2$ in $G_0$.
Let $e_1,\ldots,e_t$ be the edges of $G_0-\{e_0\}$ containing $\bar{w}$, where $t \ge 1$.
Moving those edges $e_1,\ldots,e_t$ from $\bar{w}$ to $u$, we arrive at a hypergraph $\bar{G}_0$ isomorphic to $G_0$ and $\bar{G}$ isomorphic to $G$.
However, as $x_u > x_{\bar{w}}$, $\rho(\bar{G})>\rho(G)$ by Theorem \ref{moveedges}; a contradiction.

{\scshape \mdseries Claim 2}:  The hypertree $T_u$ is a hyperstar with $u$ as its center.
If not, there exists a pendant edge $e$ of $T_u$ attached at a vertex $w \ne u$.
Relocating the edge $e$ from $w$ to $u$, we will get a hypergraph with a larger spectral radius. \hfill $\blacksquare$

For a unicyclic linear (or power) hypergraph $U$ with $m$ edges, if $m=3$ then $U$ is exactly a linear cycle of length $3$.
The {\it girth} of a hypergraph is the minimum length of the its cycles.
If a hypergraph contains no cycles, then its girth is defined to be infinity.
Denote by $S_{m,g}$ the unicyclic simple graph obtained from a cycle $C_g$ of length $g$ by attaching a star $K_{1,m-g}$ at some of its vertices.

\begin{coro} \label{maxlinegirth} Among all unicyclic linear (power) hypergraphs with $m >3$ edges and girth $g$, the power hypergraph $S_{m,g}^k$ is the unique maximizing hypergraph.
 \end{coro}

 \noindent{\bf Proof.} Let $G$ be a maximizing unicyclic linear hypergraph and let $C: v_1 e_1 v_2 \cdots v_g e_g v_1$ be the unique cycle of $G$.
 By Lemma \ref{max}, $G$ is obtained from $C$ by attaching a hyperstar at some vertex say $u$ of $C$, where $u$ has degree greater than one.
 Hence $G$ is surely a power hypergraph. The result now follows.\hfill $\blacksquare$

 \begin{lemma} For $g \ge 4$, $\rho(S_{m,g}^k) < \rho(S_{m,g-1}^k).$
 \end{lemma}

 \noindent{\bf Proof.} Let $C: v_1 e_1 v_2 \cdots v_g e_g v_1$ be the cycle of $S_{m,g}^k$, where $v_1$ is attached by a hyperstar.
Let $x$ be a Perron vector of $S_{m,g}^k$.
As shown in the proof of Corollary \ref{maxlinegirth}, $x_{v_1}> x_{v_t}$ for any $t=2,3, \ldots,g$.
Now moving the edge $e_{g-1}$ from $v_g$ to $v_1$, by Theorem \ref{moveedges} we will get a hypergraph (i.e. $S_{m,g-1}^k$) with a larger spectral radius.
\hfill $\blacksquare$

 \begin{coro} Among all unicyclic linear (power) hypergraphs with $m >3$ edges, $S_{m,3}^k$ is the unique maximizing hypergraph.
 \end{coro}

Finally we discuss the maximizing linear or power bicyclic hypergraphs.
A linear bicyclic hypergraph has at least $4$ edges; and if it has $4$ edges, then it is the hypergraph $G_5$ or $G_6$ in Fig. \ref{bic-line}.
A power bicyclic hypergraph has at least $5$ edges; and if it has $5$ edges, then it is the hypergraph $G_{10}$ in Fig. \ref{bic-line}.

\begin{theorem} \label{powerB}Among all power bicyclic hypergraph with $m>5$ edges, the unique maximizing hypergraph denoted by $B_m^P$ is obtained from the hypergraph $G_{10}$ in Fig. \ref{bic-line} by attaching a hyperstar $K_{1,m-5}$ with its center at a vertex of degree $3$.
\end{theorem}

\noindent{\bf Proof.} Let $G$ be a maximizing power bicyclic hypergraph.
Then by Lemma \ref{max}, $G$ is obtained from a bicyclic hypergraph $G_0$ by attaching a hyperstar with its center at some vertex $u$ of $G_0$, where $G_0$ consists of two cycles connected by a path (or sharing a common vertex), or three internal disjoint paths with common endpoints at most one of which has length one.
Let $x$ be a Perron vector of $G$.
 By Lemma \ref{max}, $x_u >x_v$ for any other vertex $v$.

{\scshape \mdseries Claim 1}: Every induced (or chordless) cycle of $G_0$ has length $3$.
 Otherwise, assume $G_0$ contains a cycle $C_t: v_1 e_1 v_2 \cdots e_t v_1$, where $t \ge 4$.
 Without loss of generality, assume $x_{v_1}=\max\{x_{v_i}: i=1,2,\ldots,t\}$.
 Moving the edge $e_2$ from $v_2$ to $v_1$, by Theorem \ref{moveedges}, we will get a bicyclic hypergraph with larger spectral radius.

{\scshape \mdseries Claim 2}: Every two cycles share a common vertex.
Otherwise, let $C$ and $\bar{C}$ be two cycles connected by a path with one endpoint $w$ in $C$ and the another endpoint $\bar{w}$ in $\bar{C}$.
Without loss of generality, assume $x_{w} \ge x_{\bar{w}}$.
Write $G$ as $G_1(\bar{w}) \ast G_2(\bar{w})$, where $G_1$ contains the vertex $w$.
Relocating $G_2$ from $\bar{w}$ to $w$, by Corollary \ref{relocate}, we will get a bicyclic hypergraph with larger spectral radius.

By the Claims 1 and 2, we get that $G_0$ is the hypergraph $G_{10}$ of Fig. \ref{bic-line} or two cycles of length $3$ sharing one common vertex.
Suppose that $G_0$ is the latter hypergraph.
 Let $w_1,w_2$ be the vertices of degree $2$ in one cycle and $\bar{w}_1,\bar{w}_2$ be the vertices of degree $2$ in the other cycle.
Without loss of generality assume that $x_{w_1}$ has the maximum value among $x_{w_i},x_{\bar{w}_i}$ for $i=1,2$.
Moving the edge connecting $\bar{w}_1$ and $\bar{w}_2$ from $\bar{w}_1$ to $w_1$, by Theorem \ref{moveedges}, we will get a bicyclic hypergraph with larger spectral radius.

So $G_0$ is exactly the hypergraph $G_{10}$, and $G$ is obtained from $G_{10}$ by attaching a hyperstar $K_{1,m-5}$ with its center at the vertex $u$ of $G_{10}$.
Note that $u$ has degree greater than one in $G_{10}$ by Lemma \ref{max}.
If $u$ is a vertex of degree $2$ in $G_{10}$, letting $\bar{u}$ be the other vertex of $G_{10}$ with degree $2$ and $w$ be a vertex of degree $3$,
moving the edge connecting $\bar{u}$ and $w$ from $w$ to $u$, by Theorem \ref{moveedges}, we will get a bicyclic hypergraph with larger spectral radius.
Hence $u$ has degree $3$ in $G_{10}$ and $G=B_m^P$. \hfill $\blacksquare$

 \begin{theorem}
 Among all the linear bicyclic hypergraph with $m \ge 5$ edges, the maximizing hypergraph is among one of the three hypergraphs:
 $B_m^L(1)$, $B_m^L(2)$ and $B_m^P$, where $B_m^L(1)$ and $B_m^L(2)$ are obtained from $G_5$ in Fig. \ref{bic-line} by attaching a hyperstar $K_{1,m-4}$
 with its center respectively at the vertex of degree $3$ and an arbitrary vertex of degree $2$, and $B_m^P$ is the hypergraph as in Theorem \ref{powerB}.
 \end{theorem}

 \noindent{\bf Proof.} Let $G$ be a maximizing linear bicyclic hypergraph, and let $x$ be a Perron vector of $G$.
 First suppose $G$ is not a power hypergraph.
 So there there exists an edge of $G$, say $e$, which contains at least three vertices say $v_1,v_2,v_3$ with degree greater than one.
 We have five cases according to the common neighbors among $v_1,v_2,v_3$:

 (1) $v_1,v_2,v_3$ have a common neighbor $w$ outside $e$ (see $G_5$ in Fig. 2.2);

 (2) $v_1$ and $v_2$ have a common neighbor $w$, $v_2$ and $v_3$ have a common neighbor $\bar{w}$, both outside $e$, but $w$ and $\bar{w}$ are contained in the same edge (see $G_6$ in Fig. 2.2);

 (3) $v_1$ and $v_2$ have a common neighbor $w$, $v_2$ and $v_3$ have a common neighbor $\bar{w}$, both outside $e$, but $w$ and $\bar{w}$ are not contained in the same edge (see $G_7$ in Fig. 2.2);

 (4) only $v_1$ and $v_2$ have a common neighbor $w$ outside $e$ (see $G_8$ in Fig. 2.2);

 (5) any two vertices of $v_1,v_2,v_3$ have no common neighbors (see $G_9$ in Fig. 2.2).

 We assert that only the Case (1) occurs and the other cases cannot happen.
 Suppose that Case (2) or Case (3) occurs.
 If $x_w \ge x_{\bar{w}}$, moving the edge connecting $v_3$ and $\bar{w}$ from $\bar{w}$ to $w$, we will get a bicyclic linear hypergraph with larger spectral radius
   by Theorem \ref{moveedges}.
 If $x_w < x_{\bar{w}}$, moving the edge connecting $v_1$ and $w$ from $w$  to $\bar{w}$, we also get a bicyclic linear hypergraph with larger spectral radius by Theorem \ref{moveedges}.

 If Case (4) or Case (5) occurs, and $x_{v_1} \ge x_{v_3}$ (or $x_{v_1} < x_{v_3}$), moving one edge containing $v_3$ except $e$ from $v_3$ to $v_1$
 (or moving the edge connecting $v_1$ and $w$ from $v_1$ to $v_3$), we will get a bicyclic linear graphs with larger spectral radius by Theorem \ref{moveedges}.

 So, by Lemma \ref{max}, $G$ is obtained from $G_5$ by attaching a hyperstar $K_{1,m-4}$ with its center to a vertex of degree greater than one,
 i.e. $G$ is $B_m^L(1)$ or $B_m^L(2)$.
  If $G$ is a power hypergraph, by Theorem \ref{powerB}, $G$ is the hypergraph $B_m^P$.
 The result follows.\hfill $\blacksquare$

\vspace{3mm}
{\bf Conjecture:}
For $m \ge 4$, $\rho(B_m^L(1))>\rho(B_m^L(2))$;
for $m \ge 5$, $\rho(B_m^L(1))>\rho(B_m^L(2))>\rho(B_m^P)$.



 \end{document}